\def\DATE{Sun, 20 June 2004}

\documentclass[12pt]{article}

\usepackage{amsmath,amssymb}

\allowdisplaybreaks
\newcommand\C{\mathbb{C}}

\newcommand\Aut{\mathop{\mathrm{Aut}}\nolimits}

\newtheorem{theorem}{Theorem}
\newtheorem{prop}[theorem]{Proposition}
\newtheorem{lemma}[theorem]{Lemma}

\newtheorem{example}[theorem]{Example}
\newtheorem{definition}[theorem]{Definition}
\newtheorem{rem}[theorem]{Remark}

\numberwithin{theorem}{section}
\numberwithin{equation}{section}

\makeatletter
\newenvironment{proof}[1][\proofname]{\par
  \normalfont
  \topsep6\p@\@plus6\p@ \trivlist
  \item[\hskip\labelsep{\bfseries
#1}\@addpunct{.}]\ignorespaces
}{%
  \endtrivlist
}
\newcommand{\proofname}{Proof}
\makeatother

\makeatletter
\def\qedsymbol{\RIfM@\bgroup\else$\bgroup\aftergroup$\fi
  \vcenter{\hrule\hbox{\vrule\@height.8em\kern.6em\vrule}\hrule}\egroup}
\def\qed{\RIfM@\else\unskip\nobreak\fi\quad\qedsymbol}
\makeatother

\newcommand\Pain[1]{{\mathrm P}_{\mathrm{#1}}}

\newcommand\PII{\Pain{II}}
\newcommand\PIII{\Pain{III}}
\newcommand\PIV{\Pain{IV}}
\newcommand\PV{\Pain{V}}
\newcommand\PVI{\Pain{VI}}
\newcommand\PJ{\Pain{J}}

\begin{document}

\title{Quantum Painlev\'e systems of type $A_l^{(1)}$}
\author{Hajime Nagoya \\
  {\normalsize Mathematical Institute, Tohoku University 
    Sendai 980-8578, Japan} \\
  {\normalsize\tt sa1m19@math.tohoku.ac.jp}}
\date{\DATE}
\maketitle


\begin{abstract}
We propose quantum Painlev\'e systems of type $A_l^{(1)}$.
 These systems, for $l=1$ and $l\ge 2$, should be regarded as
 quantizations of the second Painlev\'e equation
and the differential systems with the affine Weyl group symmetries 
of type $A_l^{(1)}$ studied
by M.~Noumi and Y.~Yamada \cite{NYhigherorder}, respectively.
These quantizations enjoy the affine Weyl group symmetries of
type $A_l^{(1)}$  as well as the Lax representations.
The quantized systems of type $A_1^{(1)}$ and type $A_l^{(1)}$ ($l= 2n$) 
can be obtained
as the continuous limits of the discrete systems constructed from
the affine Weyl group symmetries of type $A_2^{(1)}$ and $A_{l+1}^{(1)}$,
respectively.
\end{abstract}


\tableofcontents


\section{Introduction}
 The Painlev\'e equations $\PJ$ (J = I,...,VI) 
  were discovered by P.~Painlev\'e 
 and B.~Gambier in the classification of second-order nonlinear 
 ordinary differential equations
 without movable singular points in their solutions 
 \cite{Painleve}, \cite{Gambier}. 
  It is established by K.~Okamoto that the Painlev\'e
 equations $\PII$, $\PIII$, $\PIV$, $\PV$, and $\PVI$ have the affine 
 Weyl group symmetries of type $A^{(1)}_1$, $C^{(1)}_2$, 
  $A^{(1)}_2$, $A^{(1)}_3$, and $D^{(1)}_4$, respectively,
  as the group of B\"acklund transformations \cite{okamoto}. He also revealed
 the Hamiltonian structures for the Painlev\'e equations completely, namely,  
 we can consider the Painlev\'e equations as 
 Hamiltonian systems with affine Weyl group symmetries. 
 Therefore, a question naturally arises:
 Does there exist a quantization of Painlev\'e 
 equations with affine Weyl group symmetries?
  What we mean by the quantization is
  the canonical quantization, that is, a Poisson bracket will be 
 replaced with a commutator.
  The present paper aims to answer this problem affirmatively.
  
More on backgrounds. A.~P.~Veselov and A.~B.~Shabat studied the dressing chains 
\cite{Dressingchains}, 
and noticed that they can be considered as higher analogues of $\PIV$ 
and $\PV$. Using the dressing chains, V.~E.~Adler introduced the symmetric form of 
 the fourth Painlev\'e equation $\PIV$  \cite{Adler}. 
 Independently, using the invariant divisors of $\PIV$, M. Noumi and Y. Yamada 
 introduced the symmetric form of 
 the fourth Painlev\'e equation $\PIV$
 \cite{NYfourthokamoto}. 
 The symmetric form of $\PIV$ enables us to clarify 
 the structure of 
 the affine Weyl group symmetry. 
 Consequently, generalizing the B\"acklund transformations  
 of the symmetric form of $\PIV$, 
 M. Noumi and Y. Yamada succeeded in
 constructing a new representation of the Coxeter group associated 
 with an arbitrary generalized Cartan matrix 
  \cite{NYAffineWeyl}.  Moreover, they proposed  
 nonlinear ordinary differential systems 
 with the affine Weyl group symmetry of type  $A_l^{(1)}$ ($l\ge 2$)
  \cite{NYhigherorder}. 
  These differential systems are equivalent 
 to $\PIV$ and $\PV$ for $l=2$ and $l=3$, respectively, and they have 
   polynomial Hamiltonians \cite{NYhigherorder} and  Lax 
 representations \cite{ny1}.

  Let us formulate our main theorem.
For $l=1$, let $\mathcal{K}_1$ be the skew field
over $\mathbb{C}$ with generators $f_0$, $f_1$, $f_2$, $\alpha_0$, $\alpha_1$, $h$
and the fundamental relations
\begin{align}
&[f_1,f_0] = 2hf_2,\quad [f_0,f_2]=[f_2,f_1]=h,
\\
& 
[f_i,\alpha_j]=[h,\alpha_j]=0,
\\
&[f_i,h]=0,
\end{align} 
 and for each $l=2,3,\ldots$, let $\mathcal{K}_l$ be the skew field
   over $\mathbb{C}$ with generators  $f_i$, $\alpha_i$ ($0\le i\le l$), $h$,  
  and the fundamental relations 
    \begin{align}
&[f_i,f_{i+1}]=h,\quad [f_i,f_j]=0\quad (j\neq i\pm 1),\quad[f_i,\alpha_j]=0,
\\
& [\alpha_i,\alpha_j]=[h,f_j]=[h,\alpha_j]=0 \quad (0\le i,j\le l),
    \end{align}
where the indices $0,1,\ldots,l$ are understood as elements 
of $\mathbb{Z}/(l+1)\mathbb{Z}$. 
Then, the quantization problem of the Painlev\'e 
systems of type $A_l^{(1)}$ can be solved as follows:
  
  \begin{theorem}\label{nagoya01} 
 Let the $\mathbb{C}$-derivation $\partial$ of $\mathcal{K}_l$ be defined 
 as in Definition \ref{def:hami}. The generators of $\mathcal{K}_l$ 
 satisfy the following relations:
   
   (0) For $l=1$,
   \begin{align}
   &\partial f_0=f_0f_2+f_2f_0+\alpha_0,\quad 
   \partial f_1=-f_1f_2-f_2f_1+\alpha_1,\quad\partial f_2=f_1-f_0,
   \label{0Nl=1}
   \\
   &\partial\alpha_i=0\quad (i=0,1),\quad \partial h=0.
   \label{0Nl=1a}
   \end{align}
   
  (1)  For $l=2n$ ($n\ge 1$),
    \begin{align}
     &\partial f_i = f_i\left(\sum_{1\le r \le n}f_{i+2r-1}\right)
     -\left(\sum_{1\le r \le n}f_{i+2r}\right)f_i+\alpha_i,\label{0Nl=2n}
     \\
     &\partial\alpha_i=0\quad (0\le i \le l),\quad\partial h=0.\label{0Nl=2na}
      \end{align}
    
    (2) For $l=2n+1$ ($n\ge 1$), 
    \begin{align}
\partial f_i=&f_i\left(\sum_{1\le r\le s\le n}f_{i+2r-1}f_{i+2s}\right)
-\left(\sum_{1\le r\le s\le n}f_{i+2r}f_{i+2s+1}\right)f_i\nonumber
\\
&+\left(\frac{k}{2}-\sum_{1\le r\le n}\alpha_{i+2r}\right)f_i+
\alpha_i\sum_{1\le r\le n}f_{i+2r},\label{0Nl=2n+1}
\\
&\partial\alpha_i=0 \quad (0\le i \le l),\quad\partial h=0,\label{0Nl=2n+1a}
    \end{align}
where $k=\alpha_0+\cdots+\alpha_l$.
\end{theorem}

\begin{theorem}\label{nagoya02}   
    The $\C$-derivation $\partial$ commutes with the action of 
    the extended 
   affine Weyl group 
   $\widetilde{W}=\langle s_0,\ldots, s_l,\pi \rangle$ of type $A_l^{(1)}$ 
   defined 
   in Proposition \ref{aff1}.
\end{theorem}
   
In the classical case $h=0$, 
when $l=1$, the quantum system (\ref{0Nl=1}), (\ref{0Nl=1a}) is 
nothing but the classical second Painlev\'e equation, and 
when $l\ge 2$,  the quantum systems (\ref{0Nl=2n}), (\ref{0Nl=2na}) and  
(\ref{0Nl=2n+1}), (\ref{0Nl=2n+1a})
are nothing but the classical systems proposed by M. Noumi and Y. Yamada.
We call these systems \textsl{the quantum Painlev\'e systems} of type $A_l^{(1)}$. 

Apart from the above mentioned works on the Painlev\'e equations,  
 the Painlev\'e equations can be formulated in
    the general theory of monodromy preserving deformation
     \cite{JimboMiwaUenoMonodromy}, 
 \cite{JimboMiwaMonodromy}. 
 As for the quantization  of monodromy preserving deformation, only the
 cases of Poincar\'e rank $0$  (namely, the regular singular case) and 
  Poincar\'e rank $1$ at the infinity are known. Let us briefly mention the works 
 which are relevant to these cases. 
  
 The Schlesinger equations can be viewed as deformation equations that preserve the 
 monodromy of the rational connection 
$\partial/\partial z-\sum_{i=1}^nA_i/(z-z_i)$ 
($A_i\in M_{l+1,l+1}(\C)$) 
with regular singularities.
 N. Reshetikhin introduced the generalized Knizhnik-Zamolodchikov 
 equations and noticed that the original system of 
 the Knizhnik-Zamolodchikov equations is a quantization 
 of Schlesinger equations \cite{ReshetikhinKnizhik} (see also  
  \cite{HarnadQuantum}). 
 In the case of Poincar\'e rank $1$ at the infinity 
 where the rational connection is  
$\partial/\partial z
-\left[A+\sum_{i=1}^nB_i/(z-z_i)\right]$ 
($A, B_i\in M_{l+1,l+1}(\C)$),
   a quantization is constructed 
  (\cite{BabujianKitaevGeneralized}, \cite{FelderDifferential}).

  The classical differential systems for 
  (\ref{0Nl=1}), (\ref{0Nl=1a}) 
  (\ref{0Nl=2n}), (\ref{0Nl=2na}) 
  and (\ref{0Nl=2n+1}), (\ref{0Nl=2n+1a}) 
  describe monodromy preserving deformations of rational connections with irregular 
   singularity of Poincar\'e rank $3$ ($l=1$), $2$ ($l\ge 2$) at $z=\infty$. 
 In Proposition \ref{prop:lax}, we establish that the quantum Painlev\'e systems
 of type $A_l^{(1)}$ have the Lax representation.

     This paper is organized as follows. In Section \ref{secqp}, we  
 define the Hamiltonian of the quantum Painlev\'e systems of type $A_l^{(1)}$.
 We will also redefine 
 the quantum Painlev\'e systems of type $A_l^{(1)}$ in terms of the Hamiltonian and    
 establish the affine Weyl group symmetry. Moreover, we introduce a quantum canonical 
 coordinate and rewrite the quantum Painlev\'e systems of type $A_l^{(1)}$ into 
 the Heisenberg equations and show that the quantum Painlev\'e systems of type 
 $A_l^{(1)}$ have Lax representations. 
 In Section \ref{seccon}, we construct 
 a discrete system from the action of the extended 
 affine Weyl group of type $A_l^{(1)}$ which is defined in Subsection \ref{secaff}, 
 and take the continuous limit of the discrete system for $l=2, 2n+1$. 
 When $l=2$, we obtain the quantum second Painlev\'e equation as the continuous 
 limit, and when $l=2n+1$, we obtain the quantum Painlev\'e systems of type 
 $A_{2n}^{(1)}$.  
  See Remark \ref{rem:2n+1} for the discrete system whose continuous limit is  
  the quantum Painlev\'e systems of type 
 $A_{2n+1}^{(1)}$.

\section{Quantum Painlev\'e systems of type $A_l^{(1)}$}\label{secqp}
    For $l=1$, we can define the skew field $\mathcal{K}_1$ 
    over $\mathbb{C}$ with the generators
    \begin{equation}
f_0, f_1, f_2,\alpha_0, \alpha_1, h,
\end{equation}
and the following relations
\begin{align}
&[f_1,f_0] = 2hf_2,\quad [f_0,f_2]=[f_2,f_1]=h,
\\
& 
[f_i,\alpha_j]=[h,\alpha_j]=0,
\\
&[f_i,h]=0,
\end{align} 
       and for each $l=2,3,\ldots$, we can define the skew field $\mathcal{K}_l$ 
   over $\mathbb{C}$ with the generators 
      \begin{equation}
f_i, \alpha_i\quad (0\le i\le l), h  
 \end{equation}
  and the following relations 
    \begin{align}
    &[f_i,f_{i+1}]=h,\quad [f_i,f_j]=0\quad (j\neq i\pm 1),\quad[f_i,\alpha_j]=0,
\\
& [\alpha_i,\alpha_j]=[h,f_j]=[h,\alpha_j]=0 \quad (0\le i,j\le l),
    \end{align}
where the indices $0,1,\ldots,l$ are understood as elements 
of $\mathbb{Z}/(l+1)\mathbb{Z}$. 
  We will also identify the generators $\alpha_0,\ldots,\alpha_l$ with 
  the simple roots of the affine root system of type $A_l^{(1)}$.
  The associative algebra defined with the above relations   
    is an Ore domain, and $\mathcal{K}_l$ 
    is its quotient skew field (see, for example, \cite{bjork} Chapter 1, Section 8).

     \subsection{Hamiltonian}

Let us begin with the Hamiltonian which 
reproduces the quantum Painlev\'e systems of type $A_l^{(1)}$. 
In the classical case $h=0$, this Hamiltonian is nothing but the polynomial 
Hamiltonian for the classical Painlev\'e system of type 
$A_l^{(1)}$. Accordingly, we follow the notation of 
\cite{NYhigherorder}. 

For each $i=1,\ldots,l$, we denote by $\varpi_i$ the $i$-th fundamental weight of 
the finite root system of type $A_l$,
 \begin{align}
\varpi_i&= \frac{1}{l+1}\{(l+1-i)\sum_{r=1}^ir\alpha_r+i\sum_{r=i+1}^l(l+1-r)\alpha_r\}
\nonumber
\\
&=\sum_{r=1}^l(\mathrm{min}\{i,r\}-\frac{ir}{l+1})\alpha_r
\end{align}
and set $\varpi_0=0$. 

Put $\Gamma=\mathbb{Z}/(l+1)\mathbb{Z}$.     
For each  subset $C_{i,m}:=\{i,i+1,\ldots,i+m-1\}$ ($m\in\mathbb{Z}_{>0}$, $m\le l$)
of $\Gamma$,    
 we define $\chi(C_{i,m})$ by  
\begin{equation}
\chi(C_{i,m}):=\varpi_i-\varpi_{i+1}+\cdots+(-1)^{m-1}\varpi_{i+m-1}.
\end{equation}
For each proper subset $C=\coprod_i C_{i,m_i}$ (disjoint union) of $\Gamma$, 
 we define $\chi(C)$ by
\begin{equation}
\chi(C):=\sum_i\chi(C_{i,m_i}),
\end{equation}
where we assume that the intersection of $C_{i,m_i+1}$ and $C_{j,m_j}$
is empty  
 for $i\neq j$. Then, we call each $C_{i,m_i}$ a connected component of $C$ with
 length $m_i$.
  
For each $d=1,\ldots,l+1$, let
$S_d$ be the set of the subset $K\subset\Gamma$, such that   
 $|K|=d$, and the length of each connected component of $\Gamma\backslash K$
 is even.

For $C_{i,m}$, we set  
\begin{equation}
f_{C_{i,m}}=f_if_{i+1}\cdots f_{i+m-1}.
\end{equation}
 For each $K=\coprod_i C_{i,m_i}\in S_d$ ($d=1,\ldots,l$), 
we define $f_K$ by 
\begin{equation}
f_K=\prod_i f_{C_{i,m_i}},
\end{equation}
where $C_{i,m_i}$ is a connected component of $K$.
Note that we do not define $f_K$ for $K\in S_{l+1}$.

\begin{definition}
We define the Hamiltonian $H_0$ for the quantum Painlev\'e 
systems of type $A_l^{(1)}$ (\ref{0Nl=2n}), (\ref{0Nl=2n+1}) as follows:  

 (1) For even $l$: 
  \begin{equation}\label{hal=2n}
H_0=  \left\{
\begin{array}{lc}
f_0f_1f_2+hf_1+\sum_{K\in S_1}\chi(\Gamma \backslash K)f_K
&(l=2),
\\[3mm]
\sum_{K\in S_3}f_K+\sum_{K\in S_1}\chi(\Gamma \backslash K)f_K&
(l=2n, n\ge 2).
\end{array}\right.
\end{equation} 

(2) For odd $l$: 
\begin{equation}\label{hal=2n+1}
H_0=\left\{
\begin{array}{l}
\frac{1}{2}(f_0f_1+f_1f_0)+\alpha_1f_2
\\
(l=1)
\\
f_0f_1f_2f_3+hf_1f_2
+\sum_{K\in S_2}\chi(\Gamma \backslash K)f_K+\left(\sum_{i=1}^3(-1)^{i-1}\varpi_i\right)^2
\\(l=3),
\\[3mm]
\sum_{K\in S_4}f_K+\sum_{K\in S_2}\chi(\Gamma \backslash K)f_K+
\left(\sum_{i=1}^l(-1)^{i-1}\varpi_i\right)^2
\\(l=2n+1, n\ge 2).
\end{array}\right.
\end{equation}
\end{definition}

 Constant terms $\left(\sum_{i=1}^l(-1)^{i-1}\varpi_i\right)^2$ 
 in (\ref{hal=2n+1}) are so chosen that $H_0$ has 
 the invariance  under the affine
 Weyl group action (\ref{saction}) (see Proposition \ref{prop:hamiltonian0}).
 
\begin{example}
The explicit forms $H_0$ for $l=2,3,4,5$ are as follows:

For $l=2$:
\begin{align*}
H_0= & \ f_0f_1f_2+hf_1+\frac{1}{3}(\alpha_1-\alpha_2)f_0
+\frac{1}{3}(\alpha_1+2\alpha_2)f_1-\frac{1}{3}(2\alpha_1+\alpha_2)f_2.
\end{align*} 

For $l=3$:
\begin{align*}
H_0=&\ f_0f_1f_2f_3+hf_1f_2
+\frac{1}{4}(\alpha_1+2\alpha_2-\alpha_3)f_0f_1
+\frac{1}{4}(\alpha_1+2\alpha_2+3\alpha_3)f_1f_2
\\&
-\frac{1}{4}(3\alpha_1+2\alpha_2+\alpha_3)f_2f_3
+\frac{1}{4}(\alpha_1-2\alpha_2-\alpha_3)f_3f_0
+\frac{1}{4}(\alpha_1+\alpha_3)^2.
\end{align*}

For $l=4$:
\begin{align*}
H_0=& \ f_0f_1f_2+f_1f_2f_3+f_2f_3f_4+f_3f_4f_0+f_4f_0f_1
+\frac{1}{5}(2\alpha_1-\alpha_2+\alpha_3-2\alpha_4)f_0
\\
&+\frac{1}{5}(2\alpha_1+4\alpha_2+\alpha_3+3\alpha_4)f_1
-\frac{1}{5}(3\alpha_1+\alpha_2-\alpha_3+2\alpha_4)f_2
\\
&+\frac{1}{5}(2\alpha_1-\alpha_2+\alpha_3+3\alpha_4)f_3
-\frac{1}{5}(3\alpha_1+\alpha_2+4\alpha_3+2\alpha_4)f_4.
\end{align*}

For $l=5$:
\begin{align*}
H_0=&\ f_0f_1f_2f_3+f_1f_2f_3f_4+f_2f_3f_4f_5+f_3f_4f_5f_0+f_4f_5f_0f_1+f_5f_0f_1f_2
\\
&+\frac{1}{3}(\alpha_1+2\alpha_2+\alpha_4-\alpha_5)f_0f_1
+\frac{1}{3}(\alpha_1+2\alpha_2+3\alpha_3+\alpha_4+2\alpha_5)f_1f_2
\\&
-\frac{1}{3}(2\alpha_1+\alpha_2-\alpha_4+\alpha_5)f_2f_3
+\frac{1}{3}(\alpha_1-\alpha_2+\alpha_4+2\alpha_5)f_3f_4
\\&
-\frac{1}{3}(2\alpha_1+\alpha_2+3\alpha_3+2\alpha_4+\alpha_5)f_4f_5
+\frac{1}{3}(\alpha_1-\alpha_2-2\alpha_4-\alpha_5)f_5f_0
\\
&
+\frac{1}{3}(\alpha_1-\alpha_2+\alpha_4-\alpha_5)f_0f_3
+\frac{1}{3}(\alpha_1+2\alpha_2+\alpha_4+2\alpha_5)f_1f_4
\\&
-\frac{1}{3}(2\alpha_1+\alpha_2+2\alpha_4+\alpha_5)f_2f_5
+\frac{1}{4}(\alpha_1+\alpha_3+\alpha_5)^2.
\end{align*}
\end{example}

\begin{definition}\label{def:hami}
 Put $k=\alpha_0+\cdots+\alpha_l$. 

(1) For $l=1,2n$, we define the $\C$-derivation $\partial$ of $\mathcal{K}_l$ as follows:
\begin{align}
&\partial f_i=\frac{1}{h}[H_0,f_i]+\delta_{i,0}k,\label{hamil=2n}
\\ 
&\partial \alpha_i=\frac{1}{h}[H_0,\alpha_i],\quad
(0\le i \le l),\quad \partial h=\frac{1}{h}[H_0,h].\label{hamil=2na}
\end{align}
 
 (2) For $l=2n+1$, we define the $\C$-derivation $\partial$ of $\mathcal{K}_l$ as follows:
 \begin{align}
 &\partial f_i=\frac{1}{h}[H_0,f_i]-(-1)^i\frac{k}{2}f_i
+\delta_{i,0}kx_0,\label{hamil=2n+1}
\\
& \partial \alpha_i=\frac{1}{h}[H_0,\alpha_i]
\quad(0\le i \le l),\quad
\partial h=\frac{1}{h}[H_0,h],
\label{hamil=2n+1a}
\end{align}
where $x_0=f_0+f_2+\cdots+f_{l-1}$.
\end{definition}

We can check that $\partial$ is the $\C$-derivation of $\mathcal{K}_l$  
and that $x_0$ is a central element of $\mathcal{K}_l$ 
from the definition of $\mathcal{K}_l$.

\begin{proof}[Proof of Theorem \ref{nagoya01}]
Case $l=1$: It is straightforward to show that the right hand side of (\ref{hamil=2n})
equals to the right hand side of (\ref{0Nl=1}) from the definition (\ref{hal=2n+1}).

Case $l\ge 2$: 
For each $i=0,\ldots,l$, we can define the $\mathbb{C}$-derivation $\partial_i$ 
of $\mathcal{K}_l$ by 
\begin{equation}
\partial_i f_j=\delta_{ij},\quad \partial_i \alpha_j=0,\quad \partial_i h=0.
\end{equation}
Then, for any $\varphi\in \mathcal{K}_l$ we have 
\begin{equation}
[f_i,\varphi]=h(\partial_{i+1}-\partial_{i-1})\varphi.
\end{equation}
Using the derivation $\partial_i$, we compute $[H_0,f_i]$. 
 For the cases where $l=2,3$, we can easily calculate $[H_0,f_i]$ 
from the definition (\ref{hal=2n}), (\ref{hal=2n+1}), 
 and
then from (\ref{hamil=2n}), and (\ref{hamil=2n+1}), we obtain 
the relations  (\ref{0Nl=2n}), (\ref{0Nl=2n+1}), respectively. 
We now consider the case where $l=2n$ ($n\ge 2$). 
 For $C\subset \Gamma$, we define $S_d(C)$ by
\begin{equation}
S_d(C):=\{K\subset C\ |\ |K|=d,\ C\backslash K=\coprod_{m_i:even}C_{i,m_i}\}. 
\end{equation}
From the definition (\ref{hal=2n}), 
 we compute $\frac{1}{h}[H_0,f_i]$ as follows:
\begin{align*}
\frac{1}{h}[H_0,f_i]
&=(\partial_{i-1}-\partial_{i+1})H_0
\\
&=\sum_{K\in S_2(\Gamma\backslash\{i-1\})} f_K
-\sum_{K\in S_2(\Gamma\backslash\{i+1\})} f_K
+\chi(\Gamma\backslash\{i-1\})
-\chi(\Gamma\backslash\{i+1\})
\\
&=f_i\left(\sum_{r=1}^nf_{i+2r-1}\right)-\left(\sum_{r=1}^nf_{i+2r}\right)f_i
+\alpha_i-\delta_{i,0}k.
\end{align*}
Thus, we obtain the relations (\ref{0Nl=2n}).  We can prove the case where $l=2n+1$ 
($n\ge 2$) in a similar way. \qed
\end{proof}  

\begin{rem}
The relations in (\ref{0Nl=1}) and 
the $l+1$ relations in (\ref{0Nl=2n}) (as well as in (\ref{0Nl=2n+1})) are 
dependent among themselves. Namely, 
it holds that 
\begin{align}
&\partial(f_0+f_1+f_2^2)=k\quad (l=1),
\\
&\partial\left(\sum_{r=0}^l f_r\right)=k \quad(l=2n),\label{partiale}
\\
&\partial\left(\sum_{r=0}^n f_{2r}\right)=\frac{k}{2}\sum_{r=0}^nf_{2r},\quad
\partial\left(\sum_{r=0}^n f_{2r+1}\right)=\frac{k}{2}\sum_{r=0}^nf_{2r+1}
\quad (l=2n+1).\label{partialo}
\end{align}
\end{rem}

\subsection{Affine Weyl group symmetry}\label{secaff}
We will establish the affine Weyl group symmetry for the quantum Painlev\'e system 
of type $A_l^{(1)}$, which generalizes the symmetry for the classical system.
     Let $A=(a_{ij})_{i,j=0}^l$ be the generalized Cartan matrix of type $A_l^{(1)}$ 
 \begin{equation}\label{cartanal1}
a_{ii}=2,\quad a_{i,i\pm 1}=-1,\quad a_{ij}=0 \ (j\neq i,i\pm 1),
\end{equation}
and let $U=(u_{ij})_{i,j=0}^l$ be the matrix defined by 
\begin{equation}\label{orientational1}
u_{i,i\pm 1}=\pm 1,\quad u_{ij}=0 \ (j\neq i\pm 1).
\end{equation}    

       \begin{prop}\label{aff1}
We can define  automorphisms $s_0,\ldots,s_l, \pi\in\Aut_{\C}\mathcal{K}_l$ as follows:
 
 (1) For $l=1$,
 \begin{align}
 &s_0(f_0)=f_0,\quad 
 s_0(f_1)=f_1-f_2\frac{\alpha_0}{f_0}-\frac{\alpha_0}{f_0}f_2-\frac{\alpha_0^2}{f_0^2},
 \quad s_0(f_2)=f_2+\frac{\alpha_0}{f_0},\nonumber
 \\
 &s_1(f_0)=f_0+f_2\frac{\alpha_1}{f_1}+\frac{\alpha_1}{f_1}f_2-\frac{\alpha_1^2}{f_1^2},
 \quad s_1(f_1)=f_1,
 \quad s_1(f_2)=f_2-\frac{\alpha_1}{f_1}, \nonumber
 \\
 &s_0(\alpha_0)=-\alpha_0,\quad 
 s_0(\alpha_1)=\alpha_1+2\alpha_0,
 \quad s_1(\alpha_0)=\alpha_0+2\alpha_1,\quad s_1(\alpha_1)=-\alpha_1,\nonumber
 \\
 &\pi(f_0)=f_1,\quad \pi(f_1)=f_0,\quad \pi(f_2)=-f_2,\quad\pi(\alpha_0)=\alpha_1,\quad 
 \pi(\alpha_1)=\alpha_0,\nonumber
 \\
 &s_0(h)=s_1(h)=\pi(h)=h.\label{saction2}
 \end{align}
 (2) For $l\ge 2$,
    \begin{align}
& s_i(f_j)=f_j+\frac{\alpha_i}{f_i}u_{ij},\quad
s_i(\alpha_j)=\alpha_j-\alpha_ia_{ij},\quad s_i(h)=h,\nonumber
\\
& \pi(f_j)=f_{j+1},\quad \pi(\alpha_j)=\alpha_{j+1},\quad \pi(h)=h,
\quad(i,j\in\mathbb{Z}/(l+1)\mathbb{Z}).\label{saction}
\end{align}
    \end{prop}

    \begin{theorem}\label{aff2}
The automorphisms $s_0,\ldots s_l,\pi$ define a representation of the extended 
affine Weyl group $\widetilde{W}=\langle s_0,\ldots s_l,\pi \rangle$ 
of type $A_l^{(1)}$. 
Namely, they satisfy the commutation relations
    \begin{align}
&s_{i}^{2} = 1, \quad (s_{i}s_{j})^{3} =1 \quad (j = i \pm 1),
\quad s_is_j=s_js_i\quad (j\neq i\pm 1),
\\& \pi^{l+1}
= 1 ,\quad \pi s_{i} = s_{i+1} \pi,
    \end{align}
    where for $l=1$ they satisfy the commutation relations except $(s_{i}s_{j})^{3} =1$.  
    \end{theorem}
    
    We can prove Proposition \ref{aff1} and Theorem \ref{aff2} by direct 
    computations.
   Note that the extended affine Weyl group 
       $\widetilde{W} = W \rtimes \{1,\pi,\ldots,\pi^{l}\}$
    is the extension  of the ordinary affine Weyl group 
    $W=\langle s_0,\ldots,s_l\rangle$ by 
    the cyclic group generated by the diagram rotation $\pi$.

   The automorphisms $s_0,\ldots,s_l$ act on the Hamiltonian $H_0$ as follows. 
    We shall deal with the action of the automorphism $\pi$ on the Hamiltonian 
    $H_0$ in the Appendix. 
    
    \begin{prop}\label{prop:hamiltonian0}
With respect to the action of $W$, the Hamiltonian $H_0$ has 
the following:

 (1) For $l=1,2n$, 
  \begin{equation}\label{haw1}
    s_i(H_0)=H_0+\delta_{i,0}k\frac{\alpha_0}{f_0}\quad (i=0,\ldots,l).
\end{equation}
  
  (2) For $l=2n+1$, 
  \begin{equation}\label{haw2}
 s_i(H_0)=H_0+\delta_{i,0}k\frac{\alpha_0}{f_0}x_0\quad (i=0,\ldots,l).
\end{equation} 
In particular, 
 the Hamiltonian $H_0$ is invariant
with respect to the action of the Weyl group $W(A_l)=\langle s_1,\ldots,s_l\rangle$.
\end{prop}

We prove this proposition through direct computations. 
For practical computations,  
 it is convenient to use the Demazure operators $\Delta_i$ ($i=0,\ldots,l$) defined by 
    \begin{equation}
    \Delta_i(\varphi)=\frac{1}{\alpha_i}(s_i(\varphi)-\varphi)\quad (\varphi\in\mathcal{K}_l).
    \end{equation}
From the above definition we can easily show that $\Delta_i$ ($i=0,\ldots,l$)
satisfy the following relations: 
    \begin{align}
    &\Delta_i(\varphi\psi) = \Delta_i(\varphi)\psi+s_i(\varphi)\Delta_i(\psi)\quad 
    (\varphi,\psi\in\mathcal{K}_l),
    \\
    &\Delta_i(\alpha_i)=-2,\quad \Delta_i(\alpha_{i\pm 1})=1,\quad \Delta_i(\alpha_j)=0\quad (j\neq i,i\pm 1),
    \\
    &\Delta_i(f_i)=0,\quad \Delta_i(f_{i\pm 1})=\pm \frac{1}{f_i},\quad \Delta_i(f_j)=0\quad (j\neq i,i\pm 1).
    \end{align}
 
\begin{proof}[Proof of Proposition \ref{prop:hamiltonian0}]
 For the case where $l=1$, we can easily calculate $s_i(H_0)$ ($i=0,1$) 
 and obtain the formulas
 (\ref{haw1}).
  
 For the cases where $l=2,3$, we can easily calculate $\Delta_i(H_0)$ and
  obtain the formulas (\ref{haw1}) and (\ref{haw2}), respectively.
In the case where $l=2n$ ($n\ge 2$), we compute $\Delta_i(H_0)$ as follows:
\begin{align*}
&\quad \Delta_i\left(\sum_{K\in S_3}f_K+\sum_{K\in S_1}\chi(\Gamma \backslash K)f_K\right)
\\
=&\quad \Delta_i(f_{i-1}f_if_{i+1})+\Delta_i\left(f_if_{i+1}\sum_{r=1}^{n-1}f_{i+2r}\right)
+\Delta_i\left(f_{i+1}\sum_{1\le r\le s\le n-1}f_{i+2r}f_{i+2s+1}\right)
\\
&\quad +\Delta_i\left(\sum_{r=1}^{n-1}f_{i+2r+1}f_{i-1}f_i\right)
+\Delta_i\left(\sum_{1\le r\le s\le n-1}f_{i+2r}f_{i+2s+1}f_{i-1}\right)
\\
&\quad +\sum_{K\in S_1}\Delta_i(\chi(\Gamma \backslash K))f_K+s_i(\chi(\Gamma\backslash\{i+1\}))\frac{1}{f_i}
-s_i(\chi(\Gamma\backslash\{i-1\}))\frac{1}{f_i}
\\
=&\quad (-f_{i+1}+f_{i-1}-\frac{\alpha_i}{f_i})+\sum_{r=1}^{n-1}f_{i+2r}
+\frac{1}{f_i}\sum_{1\le r\le s\le n-1}f_{i+2r}f_{i+2s+1}
-\sum_{r=1}^{n-1}f_{i+2r+1}
\\
&\quad -\sum_{1\le r\le s\le n-1}f_{i+2r}f_{i+2s+1}\frac{1}{f_i}
+(\chi(\Gamma\backslash\{i+1\})+\alpha_i)\frac{1}{f_i}
-(\chi(\Gamma\backslash\{i-1\})-\alpha_i)\frac{1}{f_i}
\\
&\quad +\sum_{r=1}^{2n}(-1)^{r-1}f_{i+r}
\\
=&\quad \frac{\alpha_i}{f_i}+(\varpi_{i-1}-2\varpi_i+\varpi_{i+1})\frac{1}{f_i}
=\frac{\alpha_i}{f_i}+(\alpha_i+\delta_{i,0}k)\frac{1}{f_i}=\delta_{i,0}k\frac{1}{f_i}.
\end{align*} 
Consequently, we obtain the formula (\ref{haw1}). We can prove the case where
$l=2n+1$ ($n\ge 2$) in a similar way. \qed
\end{proof}
    
    \begin{proof}[Proof of Theorem \ref{nagoya02}]
     This theorem immediately follows from the definition 
     (\ref{hamil=2n}), (\ref{hamil=2n+1}),
     and 
     Proposition \ref{prop:hamiltonian0}. \qed 
    \end{proof}
\subsection{Quantum canonical coordinate and Heisenberg equation}
In the same manner as in the classical case, 
we introduce a quantum canonical coordinate for the quantum Painlev\'e system of type $A_l^{(1)}$. 
We discuss the cases of $A_{1}^{(1)}$, $A_{2n}^{(1)}$ and $A_{2n+1}^{(1)}$ 
separately. 

\medskip

Case $A_{1}^{(1)}$: Let a new quantum coordinate system be defined by 
\begin{equation}
(q;p;x)=(f_1;f_2;f_0+f_1+f_2^2).
\end{equation}
It is easy to show that
\begin{equation}
[p,q]=h,\quad [p,x]=[q,x]=0. 
\end{equation}
$H_0$ can be rewritten as a non-commutative polynomial $H=H(q;p;x)$ in the quantum 
canonical coordinate $(q;p;x)$.
 Then, we see that the quantum Painlev\'e system of type $A_1^{(1)}$ 
 (\ref{0Nl=1}) is equivalent to the 
Heisenberg equations 
\begin{equation}
\partial q=\frac{1}{h}[H,q],\quad \partial p=\frac{1}{h}[H,p],\quad \partial x=k, 
\end{equation}

\medskip

Case $A_{2n}^{(1)}$: We define a new quantum coordinate system 
\begin{equation}
(q;p;x)=(q_1,\ldots,q_n;p_1,\ldots,q_n;x),
\end{equation}
using the following formulas 
\begin{align}
&q_i=f_{2i},\quad p_i=\sum_{r=1}^i f_{2r-1}\quad (i=1,\ldots,n),
\\
&x=f_0+f_1+\cdots+f_l.\label{x}
\end{align}
  The inverse of this coordinate transformation is given by 
\begin{align}
&f_0=x-\sum_{r=1}^nq_r-p_n,\quad f_1=p_1,\quad f_2=q_1,\label{inverse0}
\\
&f_{2i-1}=p_i-p_{i-1}, \quad f_{2i}=q_i\quad (i=2,\ldots,n).\label{inverse1}
\end{align}
It is easy to show that
\begin{equation}
[p_i,q_j]=h\delta_{ij},\quad [q_i,q_j]=[p_i,p_j]=[p_i,x]=[q_i,x]=0, 
\end{equation}
for $i,j=1,\ldots,n$. By (\ref{inverse0}) and (\ref{inverse1}), 
$H_0$ can be rewritten as a non-commutative polynomial $H=H(q;p;x)$ in the quantum 
canonical coordinate $(q;p;x)$.
 Then, we see that the quantum Painlev\'e system of type $A_l^{(1)}$ 
 (\ref{0Nl=2n}) is equivalent to the 
Heisenberg equations 
\begin{equation}
\partial q_i=\frac{1}{h}[H,q_i],\quad \partial p_i=\frac{1}{h}[H,p_i],\quad \partial x=k, 
\end{equation}
where $i=1,\ldots,n$. 

\medskip

Case $A_{2n+1}^{(1)}$:   Note that from (\ref{partialo}), 
 by putting 
 \begin{equation}
\tilde{f}_{2r}=x_0f_{2r},\quad \tilde{f}_{2r+1}=x_0^{-1}f_{2r+1}\quad (r=0,1,\ldots,n),
\end{equation}
we obtain 
\begin{equation}
\partial \tilde{f}_i =\frac{1}{h}[H_0,\tilde{f}_i]+\delta_{i,0}x_0^2 \quad 
(i=0,1,\ldots,2n+1).
\end{equation}
We introduce a new coordinate system 
\begin{equation}
(q;p;x)=(q_1,\ldots,q_n;p_1,\ldots,q_n;x_0,x_1)
\end{equation}
by the following formulas:
\begin{align}
&q_i=x_0f_{2i},\quad p_i=x_0^{-1}\sum_{r=1}^i f_{2r-1}\quad (i=1,\ldots,n),
\\
&x_0=f_0+f_2+\cdots+f_{2n},\quad x_1=f_1+f_3+\cdots+f_{2n+1}.\label{x0x1}
\end{align}
 The inverse of this coordinate transformation is given by 
\begin{align}
&f_0=x_0-x_0^{-1}\sum_{r=1}^nq_r,\quad f_1=x_0p_1,\quad f_2=x_0^{-1}q_1,\label{inverse2}
\\
&f_{2i-1}=x_0(p_i-p_{i-1}), \quad f_{2i}=x_0^{-1}q_i\quad (i=2,\ldots,n).\label{inverse3}
\end{align}
It holds that
\begin{align}
&[p_i,q_j]=h\delta_{ij},\quad [q_i,q_j]=[p_i,p_j]=0,
\\
&[p_i,x_0]=[q_i,x_0]=[p_i,x_1]=[q_i,x_1]=[x_0,x_1]=0, 
\end{align}
for $i,j=1,\ldots,n$. 
By (\ref{inverse2}) and (\ref{inverse3}), 
$H_0$ can be rewritten as a non-commutative polynomial $H=H(q;p;x)$ in the quantum 
canonical coordinate $(q;p;x)$.
 Then, we see that the quantum Painlev\'e system of type $A_l^{(1)}$ 
(\ref{0Nl=2n+1}) is equivalent to the 
Heisenberg equations 
\begin{equation}
\partial q_i=\frac{1}{h}[H,q_i],\quad \partial p_i=\frac{1}{h}[H,p_i],\quad 
\partial x_0=\frac{k}{2}x_0,\quad\partial x_1=\frac{k}{2}x_1, 
\end{equation}
where $i=1,\ldots,n$.

 The above results can be summarized as follows.

\begin{theorem}\label{the:heisenberg}
(1) The quantum Painlev\'e
system of type $A_{1}^{(1)}$ (\ref{0Nl=1}) and of type $A_{2n}^{(1)}$ (\ref{0Nl=2n}) 
is equivalent to the Heisenberg equations:
\begin{equation}
\partial q_i=\frac{1}{h}[H,q_i],\quad \partial p_i=\frac{1}{h}[H,p_i],\quad \partial x=k, 
\end{equation}
where $i=1$ ($l=1$) and $i=1,\ldots,n$ ($l\ge 2$). 

 (2) The quantum Painlev\'e system of type $A_{2n+1}^{(1)}$ (\ref{0Nl=2n+1})
 is equivalent to the Heisenberg equations:
\begin{equation}
\partial q_i=\frac{1}{h}[H,q_i],\quad \partial p_i=\frac{1}{h}[H,p_i],\quad 
\partial x_0=\frac{k}{2}x_0,\quad\partial x_1=\frac{k}{2}x_1, 
\end{equation}
where $i=1,\ldots,n$.
\end{theorem}


\subsection{Lax representation}


The classical Painlev\'e systems of type $A_l^{(1)}$ arise from 
the compatibility condition 
of the linear problem  \cite{ny1}, namely, we have the 
Lax representations for those systems. 
In this subsection, we show that the quantum Painlev\'e systems of type $A_l^{(1)}$ 
 also have Lax representations. 

Let $\mathcal{A}_1$ be the skew field over $\mathbb{C}$ 
with the generators 
\begin{equation}
f_0, f_1, f_2,u_1,u_2,\epsilon_0, \epsilon_1, h,t
\end{equation}
and the following relations
\begin{align}
&[f_1,f_0] = 2hf_2,\quad [f_0,f_2]=[f_2,f_1]=h,
\\
&[u_1,f_0]=[f_0,u_2]=[u_1,f_1]=[f_1,u_2]=h,
\\
& 
[u_1,f_2]=[u_2,f_2]=[f_i,\epsilon_j]=[h,\epsilon_j]=[t,\epsilon_j]=0,
\\
&[f_i,h]=[f_i,t]=0,
\end{align} 
and let $\mathcal{A}_l$ ($l\ge 2$) be the skew field over $\mathbb{C}$ 
with the generators 
\begin{equation}
f_i,u_i,\epsilon_i,\quad (0\le i\le l), t,h,
\end{equation}
and the following relations
\begin{align}
&[f_i,f_{i+1}] = [f_i,u_{i+1}]=[u_i,f_i]=h\label{Al1},
\\
&[f_i,f_j]=0\quad (j\neq i\pm1),\quad 
[f_i,\epsilon_j]=[u_i,\epsilon_j]=[t,\epsilon_j]=[h,\epsilon_j]=0\label{Al2},
\\
&[f_i,t]=[u_i,t]=[f_i,h]=[u_i,h]=[t,h]=0,\label{Al3}
\\
&f_i-f_{i+1} = u_i-u_{i+2},\label{Al4}
\\
&f_0+f_1+\cdots+f_l=t,\label{Al5}
\end{align}
with indices understood as elements in $\mathbb{Z}/(l+1)\mathbb{Z}$.

\begin{definition}
Let $\mathcal{A}_l[z]$ be the polynomial ring, and  
we define the elements $L,B\in M_{l+1,l+1}(\mathcal{A}_l[z])$ as follows:

(1) For $l=1$,
\begin{equation}
L=\left[
  \begin{array}{cc}
    \epsilon_1+zf_2   & f_1+z    \\
   z f_0 +z^2 & \epsilon_0-zf_2    \\
       \end{array}
\right],\quad
B=\left[
  \begin{array}{cc}
    u_1   & 1     \\
      z &  u_2  \\
  \end{array}
\right].
\end{equation}

(2) For $l\ge 2$, 
\begin{equation}
L=\left[
  \begin{array}{cccccc}
    \epsilon_1   & f_1   & 1   &    &   & \\
       & \epsilon_2   & f_2   & 1   &   & \\
       &    & \ddots   &  \ddots  & \ddots &\\
     &  &    & \ddots   &  \ddots  & 1   \\
       z&    &    &    & \epsilon_l  &f_l \\
     zf_0  &  z  &    &    &   &\epsilon_0 \\
  \end{array}
\right],\quad
B=\left[
  \begin{array}{ccccc}
    u_1   & 1   &    &    &    \\
       &  u_2  &  1  &    &    \\
       &    &  \ddots  & \ddots   &    \\
       &    &    &  u_l  &  1  \\
      z &    &    &    &  u_0  \\
  \end{array}
\right].
\end{equation}
\end{definition}

Let $\partial_z$ be the $\mathcal{A}_l$-derivation of $\mathcal{A}_l[z]$ that 
maps $z$ to 1.

\begin{prop}\label{prop:lax}
For any $\mathbb{C}$-derivation $\partial_t$ of $\mathcal{A}_l[z]$ 
that maps t to 1 and z to 0 
such that 
 \begin{equation}\label{laxformula}
[z\partial_z +L,\partial_t+B]=0, 
\end{equation}
 the following formulas hold:
 for $l=1,2n$, 
  \begin{equation}\label{lft2}
  \partial_t f_i=\partial f_i,\quad \partial_t \alpha_i=0,\quad \partial_t h=0, 
    \end{equation}
 and for $l=2n+1$, 
  \begin{equation}\label{lft3}
 \partial_t f_i= \frac{2}{t}\partial f_i,\quad \partial_t \alpha_i=0,\quad \partial_t h=0,
\end{equation}
 where $\alpha_0=1-\epsilon_1+\epsilon_0$, 
$\alpha_i=\epsilon_i-\epsilon_{i+1}$ ($1\le i\le l$), and $k=1$.
Namely,   
$\partial_t$ defines the quantum Painlev\'e system for $f_i$. 
\end{prop}

\begin{rem}
The condition (\ref{laxformula}) determines a $\C$-derivation $\partial_t$ of $\mathcal{A}_l$ 
up to the action on $u_i$. Examples of such a derivation can be constructed with the 
Hamiltonian $H_0$.
\end{rem}

\begin{proof}
When $l=1$, 
the condition (\ref{laxformula}) is equivalent to the following 
system of equations
\begin{align}
&\partial_t \epsilon_i=0,\quad\partial_tf_2=f_1-f_0,
\\
&\partial_tf_0=f_0f_2+f_2f_0+\epsilon_2-\epsilon_1+1,
\\
&\partial_tf_1=-f_1f_2-f_2f_1+\epsilon_1-\epsilon_2.
\end{align}
Hence, we have the formulas (\ref{lft2}). 

When $l\ge 2$,
the condition (\ref{laxformula}) is equivalent to the following 
system of equations
\begin{align}
&\partial_t\epsilon_i=\epsilon_i u_i-u_i\epsilon_i\label{lf1},
\\
&f_i-f_{i+1}=u_i-u_{i+2}\label{lf2},
\\
&\partial_tf_i = -u_if_i+f_iu_{i+1}+\alpha_i,\label{lf3}
\end{align}
with indices understood as elements in $\mathbb{Z}/(l+1)\mathbb{Z}$. 
From the equations  (\ref{lf1}), we have $\partial_t \alpha_i=0$, and 
 the equations (\ref{lf2}) 
are the defining relations (\ref{Al4}).  
Moreover, 
one can eliminate variables $u_i$ from the right hand side of the equation  
(\ref{lf3}) by using (\ref{Al1}), 
(\ref{Al2}), (\ref{Al3}), (\ref{Al4}), and (\ref{Al5}). 
Then, one obtains the equations (\ref{lft2}) and (\ref{lft3}). 
 We explain the procedure in detail.
Note that 
inserting (\ref{Al1}) into (\ref{lf3}), we have
\begin{equation}\label{lf22}
\partial_t f_i=f_i(-u_i+u_{i+1})-h+\alpha_i.
\end{equation}

 Case  $l=2n$: 
From (\ref{lf2}) we get 
\begin{equation}\label{lf21}
\sum_{r=1}^n(f_{i+2r-1}-f_{i+2r})
=\sum_{r=1}^n (u_{i+2r-1}-u_{i+2r+1})=-u_i+u_{i+1}.
\end{equation}
From  (\ref{lf22}) and (\ref{lf21}) we obtain 
\begin{align*}
\partial_tf_i&=f_i\sum_{r=1}^n(f_{i+2r-1}-f_{i+2r})-h+\alpha_i
\\
&=f_i\left(\sum_{r=1}^nf_{i+2r-1}\right)-\left(\sum_{r=1}f_{i+2r}\right)f_i+\alpha_i.
\end{align*}
Thus, we have (\ref{lft2}). 

 Case $l=2n+1$: 
 From (\ref{lf2}) we have
\begin{equation}
\sum_{r=0}^n(f_{2r}-f_{2r+1})=\sum_{r=0}^n(u_{2r}-u_{2r+2})=0,
\end{equation}
hence, we have 
\begin{equation}\label{lf31}
\sum_{r=0}^nf_{2r}=\sum_{r=0}^nf_{2r+1}=\frac{t}{2}.
\end{equation}
From (\ref{lf3}) we get
\begin{equation}\label{lf32}
\sum_{r=0}^n\partial_tf_{i+2r}=\sum_{r=0}^n(f_{i+2r}(-u_{i+2r}+u_{i+2r+1})
+\alpha_{i+2r}-h).
\end{equation}
For each $r=1,\ldots,n$, we define $B_r\in\mathcal{A}_l$ by the following relation:  
\begin{equation}\label{lf33}
u_{i+2r}-u_{i+2r+1}=u_i-u_{i+1}+B_r.
\end{equation}
Then, from (\ref{lf2}) we have
\begin{equation}\label{lf34}
B_r=\sum_{k=1}^rf_{i+2k-1}-\sum_{k=r+1}^{n+1}f_{i+2k-1}
-\sum_{k=1}^{r-1}f_{i+2k}+\sum_{k=r+1}^nf_{i+2k}.
\end{equation}
From (\ref{lf31}) and  (\ref{lf32}) we have
\begin{equation}\label{lf35}
-u_i+u_{i+1}=\frac{2}{t}\left\{\frac{1}{2}+\sum_{r=1}^nf_{i+2r}B_r
-\sum_{r=0}^n\alpha_{i+2r}+(n+1)h
\right\}.
\end{equation}
Inserting (\ref{lf35}) into (\ref{lf22}), we have 
\begin{equation}\label{lf36}
\frac{t}{2}\partial_tf_i=f_i\left\{\frac{1}{2}+\sum_{r=1}^nf_{i+2r}B_r
-\sum_{r=0}^n\alpha_{i+2r}+(n+1)h\right\}+\frac{t}{2}(-h+\alpha_i).
\end{equation}
Substituting (\ref{lf34}) for $B_r$, we have
\begin{align}
\frac{t}{2}\partial_tf_i&=f_i\left\{\frac{1}{2}+\sum_{r=1}^nf_{i+2r}
\left(\sum_{k=1}^rf_{i+2k-1}-\sum_{k=r+1}^{n+1}f_{i+2k-1}
-\sum_{k=1}^{r-1}f_{i+2k}+\sum_{k=r+1}^nf_{i+2k}\right)\right.\nonumber
\\
&\left.
-\sum_{r=0}^n\alpha_{i+2r}+(n+1)h\right\}
+\frac{t}{2}(-h+\alpha_i).\label{lf37}
\end{align}
Since 
\begin{equation}\label{lf38}
\sum_{r=1}^nf_{i+2r}
\sum_{k=1}^rf_{i+2k-1}=\sum_{r=1}^n\sum_{k=1}^rf_{i+2k-1}f_{i+2r}-nh,
\end{equation}
and
\begin{equation}\label{lf39}
\sum_{r=1}^nf_{i+2r}\left(-\sum_{k=1}^{r-1}f_{i+2k}+\sum_{k=r+1}^nf_{i+2k}\right)=0,
\end{equation}
we obtain
\begin{align*}
\frac{t}{2}\partial_tf_i=&f_i\left\{\frac{1}{2}+\sum_{r=1}^n\sum_{k=1}^rf_{i+2k-1}f_{i+2r}
-\sum_{r=1}^n\sum_{k=r+1}^{n+1}f_{i+2r}f_{i+2k-1}\right.
\left.
-\sum_{r=0}^n\alpha_{i+2r}+h\right\}
\\
&+\frac{t}{2}(-h+\alpha_i)
\\
=&f_i\left(\sum_{1\le k\le r\le n}f_{i+2k-1}f_{i+2r}\right)-
\left(\sum_{1\le k\le r\le n}f_{i+2r}f_{i+2k-1}\right)f_i+h\sum_{r=1}^nf_{i+2r}
\\
&+f_i\left(\frac{1}{2}-\sum_{r=0}^n\alpha_{i+2r}\right)+hf_i+\sum_{r=0}^nf_{i+2r}(-h+\alpha_i)
\\
=&f_i\left(\sum_{1\le k\le r\le n}f_{i+2k-1}f_{i+2r}\right)-
\left(\sum_{1\le k\le r\le n}f_{i+2r}f_{i+2k-1}\right)f_i
\\
&+f_i\left(\frac{1}{2}-\sum_{r=0}^n\alpha_{i+2r}\right)
+\alpha_i\sum_{r=1}^nf_{i+2r}.
\end{align*}
Thus, we have (\ref{lft3}). \qed
\end{proof}

As in the classical case \cite{noumi}, from the viewpoint of the Lax representation, 
the origin of the affine Weyl group symmetry for the quantum Painlev\'e 
systems can be explained as follows.

       Let $G_i(z),\Lambda(z)$ be matrices of 
       $M_{l+1,l+1}(\mathcal{A}_l[z,z^{-1}])$ 
       defined by 
\begin{align*}
&G_0(z)=1+\frac{\alpha_0}{f_0}z^{-1}E_{1,l+1},
\quad G_i(z)=1+\frac{\alpha_i}{f_i}E_{i+1,i} \quad (1\le i\le l),
\\
&\Lambda(z)=\sum_{i=1}^l E_{i,i+1}+zE_{1,l+1},
\end{align*}
where $E_{ij}$ is the matrix unit with $1$ at the ($i,j$) entry 
and $0$ for other entries.

We define the action of $w=s_0,\ldots,s_l,\pi$ on $\mathcal{A}_l$
as follows:
 \begin{align*}
&z\partial_z+w(L)=G_w(z)\left(z\partial_z+L\right)(G_w(z))^{-1},
\\
&\partial_t+w(B)=G_w(z)\left(\partial_t+B\right)(G_w(z))^{-1},
\end{align*}
where $G_w=G_i$ for $w=s_0,\ldots,s_l$ and $G_{\pi}=\Lambda$.
Then, from the definition, the action of 
$s_0,\ldots,s_l,\pi$ on $\mathcal{A}_l$ commutes with the derivation 
$\partial_t$ that satisfies (\ref{laxformula}). It can be seen that
the action of $s_0,\ldots,s_l,\pi$ for $f_i$, $\alpha_i$ 
is nothing but the action (\ref{saction}) of $s_0,\ldots,s_l,\pi$ on $\mathcal{K}_l$.


\section{Continuous limit}\label{seccon}

  
  \subsection{Discrete system}
  We construct a quantum discrete system with  
   affine Weyl group symmetry of type $A_l^{(1)}$ in the same way as in
    \cite{NYAffineWeyl}. We introduce the shift operators 
  $T_i$ ($1\le i\le l+1$) by 
      \begin{align*}
        &T_1=\pi s_ls_{l-1}\cdots s_1,
        \\
        &T_2=s_1\pi s_ls_{l-1}\cdots s_2,
        \\
        &\qquad\vdots
        \\
        &T_{l+1}=s_l s_{l-1}\cdots s_1 \pi.
      \end{align*}
      Then, we have the following relations: 
      \begin{align*}
    &T_iT_j=T_jT_i,\quad T_1\cdots T_{l+1}=1,
    \\
    &T_i(\alpha_{i-1})=\alpha_{i-1}+k,\quad T_i(\alpha_i)=\alpha_i-k,
    \quad T_i(\alpha_j)=\alpha_j\quad(j\neq i,i-1).
    \end{align*}
    We consider the shift operator $T_i$ as the time evolution operator. 
    Since  $T_j=\pi^{j-1}T_1\pi^{1-j}$ for $j=2,\ldots,l+1$, in the following we 
    take $T_1$ without loss of generality.
    The quantum discrete system with $T_1$ as its time evolution operator
    has the form  
    \begin{equation}\label{dissys}
    f_i[n+1]=G_i[n]\quad (0\le i\le l),
    \end{equation}
    where  $f_i[n]$ stands for $T_1^n(f_i)$, $G_i[n]$ is 
    a rational function of $f_j[n],\alpha_j[n]$ for each $i=0,1,\ldots,l$.
    
    For example, when $l=2$, the quantum discrete system is written as follows:
    \begin{equation}\label{eq:discretesystem2}
    \left\{
    \begin{array}{l}
    f_0[n+1]=
    f_1[n]+\cfrac{\alpha_0[n]}{f_0[n]}-\cfrac{\alpha_2[n]+\alpha_0[n]}
    {f_2[n]-\cfrac{\alpha_0[n]}{f_0[n]}},
    \\[3mm]
    f_1[n+1] = f_2[n]-\cfrac{\alpha_0[n]}{f_0[n]},
    \\[3mm]
    f_2[n+1] =
    f_0[n]+\cfrac{\alpha_2[n]+\alpha_0[n]}{f_2[n]-\cfrac{\alpha_0[n]}{f_0[n]}}.
    \end{array}
    \right.
    \end{equation}
    
    For each $l\ge 2$, 
    the quantum discrete system (\ref{dissys}) has the affine Weyl group symmetry of 
    type $A_{l-1}^{(1)}$, because automorphisms $s_0s_1s_0$, $s_2,\ldots,s_l$ 
    of $\mathcal{K}_l$ commute 
    with $T_1$, and $s_0s_1s_0$, $s_2,\ldots,s_l$ define a representation of 
    the affine Weyl group of type $A_{l-1}^{(1)}$.
    Therefore, if one can take an appropriate continuous limit 
    of this discrete system, 
    one would obtain a continuous system with affine 
    Weyl group symmetry of type $A_{l-1}^{(1)}$ in $\mathcal{K}_l$.
  
    Indeed, we can take an appropriate continuous limit if $l$ is either $2$ or $2n+1$ 
    ($n=1,2,\ldots$). We shall see how to take a continuous limit 
    in the next subsection.
  
  \subsection{How to take a continuous limit}
 When $l=2$, we obtain the quantum second Painlev\'e equation as the continuous limit, 
 and when $l=2n+1$, we obtain the quantum Painlev\'e system of type  
 $A_{2n}^{(1)}$ as the continuous limit as we shall see below.
 
 Informally, we consider the continuous limit as follows:
  First, we introduce the parameter $\epsilon$ called the lattice parameter, and 
 then we introduce the continuous time variable $t$ such that $t=n\epsilon$, where 
 $n$ is the discrete time variable. Second, for a function of $n$ we set 
     \begin{equation}\label{cl1}
        f[n]=y_0+y_1\epsilon + y_2\frac{\epsilon^2}{2}+\cdots,
     \end{equation} 
     where $y_i$ is a function of $t$. 
     Third,  assuming 
     \begin{equation}\label{cl2}
    f[n+1]=f[n]+\epsilon\frac{df[n]}{dt}+
    \frac{\epsilon^2}{2}\frac{d^2f[n]}{dt^2}+\cdots,
    \end{equation}
    and 
    comparing the above equation with $f[n+1]=G[n]$, where $G[n]$ is 
    a function  of $f[n]$,  we obtain the derivative $dy_i/dt$. The differential 
    equations for $dy_i/dt$ ($i=0,1,\ldots$) are  
    its continuous limit.   
    
    Now we will take the continuous limit for the case where $l=2n+1$.
    We can define the skew field 
    $\mathcal{F}_l$ over $\mathbb{C}$ with the generators 
    $\varphi_i$, $\beta_i$ ($0\le i\le l$), $t$, $h'$  
    and the following relations
    \begin{align*}
    &[\varphi_i,\varphi_{i+1}]=h',
    \quad[\varphi_i,\varphi_j]=0\quad (j\neq i\pm 1),
    \\
    &[\varphi_i,\beta_j]=0,\quad [\beta_i,\beta_j]=0,
    \\
    &[\varphi_i,t]=[\beta_i,t]=[\varphi_i,h']=[\beta_i,h']=[t,h']=0,
    \\
    &\varphi_0+\varphi_2+\cdots+\varphi_{2n}=0,\quad 
    \varphi_1+\varphi_3+\cdots+\varphi_{2n+1}=0,
    \\
    &\beta_0+\beta_1+\cdots+\beta_l=1, 
    \end{align*}
    where the indices $0,1,\ldots,l$ are understood as elements 
    of $\mathbb{Z}/(l+1)\mathbb{Z}$.  
    Also, since $f_0+f_2+\cdots+f_{2n},f_1+f_3+\cdots+f_{2n+1}
    \in\mathcal{K}_l$ 
    are central elements in $\mathcal{K}_l$ and invariants of the action of 
    the affine Weyl group $W$, we put 
    these elements as a constant in this quantum discrete system. 
    In particular, we set $f_0+f_2+\cdots+f_{2n}:=1$, 
    $f_1+f_3+\cdots+f_{2n+1}:=1$.

    \begin{lemma}Let $\mathcal{F}_l(\epsilon)$ be the quotient skew field of 
    the polynomial ring $\mathcal{F}_l[\epsilon]$ with coefficients in $\mathcal{F}_l$.  
    We can define the homomorphism $\Psi:\mathcal{K}_l\to\mathcal{F}_l(\epsilon)$ as follows:
    \begin{align*}
        &\Psi(f_0)=1+\epsilon \varphi_0,\quad \Psi(f_1)=1+\epsilon \varphi_1,
        \quad \Psi(f_i)=\epsilon\varphi_i \ (2\le i\le l),
        \\
        &\Psi(\alpha_0)=-1+\epsilon t+\epsilon^2\beta_0,\quad
        \Psi(\alpha_1)=1-\epsilon t+\epsilon^2\beta_1,
        \quad \Psi(\alpha_i)=\epsilon^2\beta_i\ (2\le i\le l),
        \\
        &\Psi(h)=\epsilon^2 h'.
    \end{align*}
    \end{lemma}

    \begin{proof}
    One can show that $\Psi$ preserves the defining relations 
    by using the definition of 
    $\mathcal{K}_l$ and $\mathcal{F}_l$. \qed
    \end{proof}
    
    We introduce the elements $\psi_i$ and $\gamma_i$ 
    of $\mathcal{F}_l$ ($0\le i \le l-1$) 
     by 
    \begin{align*}
        &\psi_0=\varphi_0+\varphi_1+t,
        \quad\psi_i=\varphi_{i+1}\quad (1\le i\le l-1),
        \\
        &\gamma_0=\beta_0+\beta_1,
        \quad
        \gamma_i=\beta_{i+1}\quad(1\le i\le l-1).
    \end{align*}
    We denote $s_0s_1s_0,s_2,\ldots,s_l$ by $r_0,r_1,\ldots,r_{l-1}$, respectively.

     We can define the action of the subgroup 
     $\widetilde{W}'=\langle T_1,r_0,\ldots,r_{l-1}\rangle$ 
    of $\widetilde{W}$ on $\mathcal{F}_l(\epsilon)$ as follows: 
        \begin{align}
    &T_1(\epsilon)=\epsilon,\quad T_1(t)=t+\epsilon,
    \quad T_1(\beta_i)=\beta_i\quad(i=0,\ldots,l),\nonumber
    \\
    & 
    T_1(\psi_i)=\frac{\Psi T_1(f_{i+1})}{\epsilon}\quad (i=1,\ldots,l-1),\label{w1}
    \end{align}
    and for each $i=0,\ldots,l-1$, 
    \begin{align}
    &r_i(\epsilon)=\epsilon,\quad r_i(t)=t,\nonumber
    \\
    &r_i(\gamma_j)=\gamma_j-\gamma_ia_{ij}\quad (j=0,\ldots,l-1),\nonumber
    \\
    &r_i(\psi_j)=\frac{\Psi r_i(f_{j+1})}{\epsilon}\quad (j=1,\ldots,l-1),\nonumber
    \\
    &r_0(\beta_0)=-\beta_1,\quad r_l(\beta_0)=\beta_0+\beta_l,\quad 
    r_i(\beta_0)=\beta_0\quad (1\le i\le l-1), \label{w2}
       \end{align}
       where $A=(a_{ij})_{i,j=0}^{l-1}$ is the generalized Cartan matrix of type 
       $A_{l-1}^{(1)}$ (\ref{cartanal1}). Then the homomorphism $\Psi$ is an 
    $\widetilde{W}'$-intertwiner.  
        
    We have the next lemma from the action of 
$\langle T_1,r_0,\ldots,r_{l-1}\rangle$ on $\mathcal{F}_l(\epsilon)$. 

\begin{lemma}\label{partial1}
Let $\mathcal{F}_l[[\epsilon]]$ be the ring of the formal power series with coefficients in 
$\mathcal{F}_l$. Then, 
we can define the action of $\widetilde{W}'$ 
on $\mathcal{F}_l[[\epsilon]]$ by the relations (\ref{w1}), (\ref{w2}). 
Furthermore, for $\phi\in\mathcal{F}_l[[\epsilon]]$,  
we have
\begin{equation}\label{time1}
T_1(\phi)=\phi+\epsilon \tilde{\phi},\quad \text{for some 
$\tilde{\phi}\in\mathcal{F}_l[[\epsilon]]$}.
\end{equation}
 \end{lemma}

\begin{proof}
 We can directly compute $r_i(\psi_j)$ and  
 obtain 
 \begin{equation}
r_i(\psi_j)=
\psi_j+\frac{\gamma_i}{\psi_i}u_{ij}+\epsilon\phi,\quad
 \text{for some $\phi\in\mathcal{F}_l[[\epsilon]]$},
\end{equation}
where $U=(u_{ij})_{i,j=0}^{l-1}$  as in (\ref{orientational1}).

For $i=1,\ldots,l-1$, one can prove inductively that 
the following formulas hold:
\begin{equation}\label{eq:conind}
\Psi \pi s_l\cdots s_{i+1}(f_i)=
\left\{
\begin{array}{lc}
1+\epsilon(\varphi_{i+1}-\varphi_{i+2}+\ldots+\varphi_l-\varphi_0-t)+\epsilon^2 \phi
&\text{for even $i$}
\\
1+\epsilon(\varphi_{i+1}-\varphi_{i+2}+\ldots-\varphi_l+\varphi_0)+\epsilon^2 \phi
&\text{for odd $i$ }
\end{array}
\right.
\end{equation}
where $\phi$ is some element of $\mathcal{F}_l[[\epsilon]]$.
Applying these formulas to 
\begin{equation}\label{eq:concom}
T_1(f_{i+1})=\pi s_l\cdots s_{i+2}(f_{i+1})+\pi s_l\cdots s_{i+1}
\left(\frac{\alpha_{i}}{f_i}\right)
\end{equation}
for $i=1,\ldots,l-1$, we obtain that  $T_1(\psi_i)=\psi_i+\epsilon\phi$ for 
some $\phi\in\mathcal{F}_l[[\epsilon]]$. As a result, we obtain the formulas (\ref{time1}).
   \qed
    \end{proof}

 From Lemma \ref{partial1}, we can define the 
 following $\C$-derivation $\partial$ 
of $\mathcal{F}_l$.
\begin{definition}\label{def:bibun}
 We define the $\C$-derivation $\partial$ of 
$\mathcal{F}_l$ by
\begin{equation}
\partial \phi=\left. \frac{T_1(\phi)-\phi}{\epsilon}\right|_{\epsilon=0},
\end{equation}
for $\phi\in\mathcal{F}_l$.
\end{definition}

From Lemma \ref{partial1}, it holds that 
$r_i\left(\epsilon \mathcal{F}_l[[\epsilon]]\right)\subset\epsilon \mathcal{F}_l[[\epsilon]]$, 
hence   
 the actions of $r_0,\ldots,r_{l-1}$ 
on $\mathcal{F}_l=
 \mathcal{F}_l[[\epsilon]]/\epsilon \mathcal{F}_l[[\epsilon]]$, 
  are induced. 
Then we have the next Theorem. 

\begin{theorem}\label{the:partial2}
(1) The $\C$-derivation $\partial$ of $\mathcal{F}_l$ acts on 
$\psi_i$ ($i=0,\ldots,l-1$) as follows:
     \begin{align}
     &\partial \psi_i = \psi_i\left(\sum_{1\le r \le n}\psi_{i+2r-1}\right)
     -\left(\sum_{1\le r \le n}\psi_{i+2r}\right)\psi_i+\gamma_i.\label{2n}
      \end{align}
    
    (2) The action of $r_0,\ldots,r_{l-1}$ on $\mathcal{F}_l$
     commutes with the derivation $\partial$.
     \end{theorem}
     
     \begin{proof}
    (1)  For $i=1,\ldots,l-1$, we define the elements $a_i,b_i\in\mathcal{F}_l$ by   
    \begin{equation}\label{eq:concom1}
\Psi \pi s_l\cdots s_{i+1}(f_i)=1+\epsilon a_i+\epsilon^2 b_i+\epsilon^3 \phi\quad
\text{for some $\phi\in\mathcal{F}_l[[\epsilon]]$}. 
    \end{equation} 
    Then, we have
    \begin{align*}
    \Psi\pi s_l\cdots s_{i+1}(f_i)
    =&\Psi\pi s_l\cdots s_{i+2}(f_i-\frac{\alpha_{i+1}}{f_{i+1}})
    \\
    =&\epsilon \varphi_{i+1}
    -\{-1+\epsilon t+\epsilon^2(\beta_{i+2}+\cdots+\beta_l+\beta_0)\}
    \\
    &\times\{1-\epsilon a_{i+1}-\epsilon^2(b_{i+1}-a_{i+1}^2)\}+\epsilon^3 \phi'
    \\
    =&1+\epsilon a_i+\epsilon^2(ta_{i+1}-(\beta_{i+2}+\cdots+\beta_l+\beta_0
    +b_{i+1}-a_{i+1}^2))+\epsilon^3 \phi'',
    \end{align*}
     where $\phi'$ and $\phi''$ are some elements of 
    $\mathcal{F}_l[[\epsilon]]$.
    Consequently, by using (\ref{eq:conind}) and (\ref{eq:concom}), we have
    \begin{align*}
    \Psi T_1(f_{i+1})
    =&\Psi(\pi s_l\cdots s_{i+2}(f_{i+1})+\pi s_l\cdots s_{i+1}(\frac{\alpha_{i}}{f_i}))
    \\
    =&1+\epsilon a_{i+1}+\epsilon^2 b_{i+1}+\{-1+\epsilon t
    +\epsilon^2(\beta_{i+1}+\cdots+\beta_l+\beta_0)\}
        \\
    &\times\{1-\epsilon a_i+\epsilon^2(\beta_{i+2}+\cdots+\beta_l+\beta_0
    +b_{i+1}-a_{i+1}^2-ta_{i+1}+a_i^2)\}+\epsilon^3 \phi'''
    \\
    =&\epsilon(a_{i+1}+a_i+t)
    +\epsilon^2(\beta_{i+1}+a_{i+1}^2-a_i^2+t(a_{i+1}-a_i))+\epsilon^3 \phi'''',
    \end{align*}
    where  $\phi'''$ and $\phi''''$ are some elements of 
    $\mathcal{F}_l[[\epsilon]]$. From (\ref{eq:conind}), for even $i$,  we have 
    \begin{align*}
    \partial \psi_i=&\beta_{i+1}+a_{i+1}^2-a_i^2+t(a_{i+1}-a_i)
    \\
    =&\gamma_i+(\varphi_{i+2}-\varphi_{i+3}+\ldots-\varphi_l+\varphi_0)^2
    -(\varphi_{i+1}-\varphi_{i+2}+\ldots+\varphi_l-\varphi_0-t)^2
    \\
    &+t(t-\varphi_{i+1}-2(-\varphi_{i+2}+\ldots+\varphi_l-\varphi_0))
    \\
    =&\gamma_i
    -\varphi_{i+1}(-\varphi_{i+2}+\ldots+\varphi_l-\varphi_0)+t\varphi_{i+1}
    \\
    &-(-\varphi_{i+2}+\ldots+\varphi_l-\varphi_0)\varphi_{i+1}-\varphi_{i+1}^2
    \\
    =&\gamma_i+t\psi_i-\psi_i\sum_{r=1}^n\psi_{i+2r}
    -\sum_{r=1}^n\psi_{i+2r}\psi_i-\psi_i^2
    \\
    =&\gamma_i+\psi_i\sum_{r=1}^n\psi_{i+2r-1}-\sum_{r=1}^n\psi_{i+2r}\psi_i.
        \end{align*}
    For odd $i$, one can compute it in a similar way. Hence, the formulas (\ref{2n}) 
    hold. 
    
             (2) The action of $r_0,\ldots,r_{l-1}$ commutes with the action of $T_1$. 
     Hence,  from the definition 
     of the derivation $\partial$, the action of $r_0,\ldots,r_{l-1}$ 
     commutes with $\partial$. 
     \qed
     \end{proof}

Hence, the derivation $\partial$ of $\mathcal{F}_{2n}$ defines the quantum Painlev\'e system of 
type $A_{2n}^{(1)}$. However, it is not clear that there is some connection between 
the derivation $\partial$ of $\mathcal{K}_{2n}$
 defined by the Hamiltonian $H_0$ in Section \ref{secqp} and 
the derivation $\partial$ of $\mathcal{F}_{2n}$ 
 defined by the time evolution $T_1$ in this section.

\begin{rem}\label{rem:2n+1}
 In the classical case, a continuous 
limit of the discrete system constructed from a representation of 
$W(A_2^{(1)})\times W(A_l^{(1)})$ ($l\ge 2$) in \cite{kny2} 
is the classical Painlev\'e system of type $A_l^{(1)}$. 
If one can quantize the representation of $W(A_2^{(1)})\times W(A_l^{(1)})$, then
one would obtain the quantum discrete system whose continuous limit is 
the quantum Painlev\'e system of type $A_l^{(1)}$. 
\end{rem}

      An appropriate continuous limit of the quantum discrete system 
      for $l=2$ (\ref{eq:discretesystem2}) is the quantum second Painlev\'e equation. 
      We can take the limit in the similar way as in the case where
      $l=2n+1$, though how to set  $y_i$ in (\ref{cl1}) is not exactly the same. 
      Let $\mathcal{F}_2$ be the skew field over $\C$ 
     with the generators $\psi,\varphi_0,\varphi_1$, $\beta_0,\beta_1,\beta_2$, $t,h'$ and 
     the following relations: 
     \begin{align*}
        &[\psi,\varphi_i]=\frac{h'}{2}\quad (i=0,1),\quad[\varphi_0,\varphi_1]=0,
        \\
        &\text{$\beta_i,t,h'$ are central},
        \\
        &\beta_0+\beta_1+\beta_2=1.
    \end{align*}
      
      Since $f_0+f_1+f_2$ is central in $\mathcal{K}_2$ and 
      the invariant of the action of $W$, we put $f_0+f_1+f_2=2$. 
      
      \begin{lemma}
      We can define the homomorphism $\Psi:\mathcal{K}_2\to\mathcal{F}_2(\epsilon)$ as follows:
    \begin{align*}
        &\Psi(f_0)=1+\epsilon \psi+\epsilon^2 \varphi_0,\quad 
        \Psi(f_1)=1-\epsilon\psi+\epsilon^2 \varphi_1,
        \quad \Psi(f_2)=-\epsilon^2(\varphi_0+\varphi_1),
        \\
        &\Psi(\alpha_0)=-1+\epsilon^2 t+\epsilon^3\beta_0,\quad
        \Psi(\alpha_1)=1-\epsilon^2 t+\epsilon^3\beta_1,
        \quad \Psi(\alpha_2)=\epsilon^3\beta_2,
        \\
        &\Psi(h)=\epsilon^3 h'.
    \end{align*}
\end{lemma}
       
    \begin{proof}
    One can show that $\Psi$ preserves the defining relations
    by using the definition of 
    $\mathcal{K}_2$ and $\mathcal{F}_2$. \qed
    \end{proof}
    
      We can define the action of the subgroup 
     $\widetilde{W}'=\langle T_1,r_0,r_1\rangle$ 
    of $\widetilde{W}$ on $\mathcal{F}_2(\epsilon)$ as follows: 
        \begin{align}
    &T_1(\epsilon)=\epsilon,\quad T_1(t)=t+\epsilon,
    \quad T_1(\beta_i)=\beta_i\quad(i=0,1,2),\nonumber
    \\
    & 
    T_1(\psi)=\psi+\epsilon(2(\varphi_0+\varphi_1)-\psi^2+t),\nonumber
    \\
    &T_1(\varphi_0)=\frac{\Psi T_1(f_0)-1-\epsilon T_1(\psi)}{\epsilon^2},
    \quad T_1(\varphi_1)=\frac{\Psi T_1(f_1)-1+\epsilon T_1(\psi)}{\epsilon^2},\nonumber
       \end{align}
    and for each $i=0,1$, 
    \begin{align}
    &r_i(\epsilon)=\epsilon,\quad r_i(t)=t,\nonumber
    \\
    &r_0(\psi)=\psi-\frac{\beta_0+\beta_1}{\varphi_0+\varphi_1+t-\psi^2},
\quad r_1(\psi)=\psi-\frac{\beta_2}{\varphi_0+\varphi_1},\nonumber
    \\
    &r_i(\varphi_0)=\frac{\Psi r_i(f_0)-1-\epsilon r_i(\psi)}{\epsilon^2},
    \quad r_i(\varphi_1)=\frac{\Psi r_i(f_1)-1+\epsilon r_i(\psi)}{\epsilon^2},\nonumber
    \\
    &r_0(\beta_0)=-\beta_1,\quad r_0(\beta_1)=-\beta_0,\quad r_0(\beta_2)=2-\beta_2,\nonumber
    \\ 
    &r_1(\beta_0)=\beta_0+\beta_2,\quad r_1(\beta_1)=\beta_1+\beta_2, 
    \quad r_1(\beta_2)=-\beta_2.\nonumber
       \end{align}
       Then, the homomorphism $\Psi$ is an 
    $\widetilde{W}'$-intertwiner, and 
we can define the action of $\widetilde{W}'$ 
on $\mathcal{F}_2[[\epsilon]]$ by the above relations. 
Furthermore, we have
\begin{align}
T_1(\varphi_0)=&\varphi_0+\epsilon(\psi\varphi_1+\varphi_1\psi+\psi^3
-t\psi+\beta_0-\beta_2)+\epsilon^2\phi,
\\ 
T_1(\varphi_1)=&\varphi_1+\epsilon(\psi\varphi_0+\varphi_0\psi-\psi^3
+t\psi-\beta_0)+\epsilon^2\phi,
\\
r_0(\varphi_0)=&\varphi_0-\frac{1}{\varphi_0+\varphi_1+t-\psi^2}\psi(\beta_0+\beta_1)
\nonumber
\\
&+\frac{1}{\varphi_0+\varphi_1+t-\psi^2}(\beta_0+\psi\varphi_1-\varphi_0\psi)
\frac{\beta_0+\beta_1}{\varphi_0+\varphi_1+t-\psi^2}+\epsilon\phi',
\\
r_0(\varphi_1)=&\varphi_1-\frac{\beta_0+\beta_1}{\varphi_0+\varphi_1+t-\psi^2}\psi
\nonumber
\\
&+\frac{1}{\varphi_0+\varphi_1+t-\psi^2}(\beta_1+\psi\varphi_0-\varphi_1\psi)
\frac{\beta_0+\beta_1}{\varphi_0+\varphi_1+t-\psi^2}+\epsilon\phi'',
\end{align}
where $\phi$, $\phi'$, and $\phi''$ are some elements of $\mathcal{F}_2[[\epsilon]]$. 
Hence, we can define 
the $\C$-derivation $\partial$ of $\mathcal{F}_2$ by 
\begin{equation}
\partial \phi=\left. \frac{T_1(\phi)-\phi}{\epsilon}\right|_{\epsilon=0},
\end{equation}
for $\phi\in\mathcal{F}_2$, and on $\mathcal{F}_2$, the actions 
of $r_0$, $r_1$ are induced. Then we have 
the next theorem. 

\begin{theorem}
(1) The $\C$-derivation $\partial$ of $\mathcal{F}_2$ acts on $\psi$, $\varphi_0$ and $\varphi_1$ as follows:
\begin{align*}
&\partial \psi=2(\varphi_0+\varphi_1)-\psi^2+t,
\\
&\partial \varphi_0=\psi\varphi_1+\varphi_1\psi+\psi^3
-t\psi+\beta_0-\beta_2,\quad \partial \varphi_1=\psi\varphi_0+\varphi_0\psi-\psi^3
+t\psi-\beta_0.
\end{align*}
Therefore, we have 
\begin{equation}\label{qp2}
\partial^2 \psi=2\psi^3-2t\psi-2\beta_2+1.
\end{equation}

(2) The actions of $r_0$, $r_1$ on $\mathcal{F}_2$ commute with the derivation $\partial$.
\end{theorem}   
       
       In the classical case $h'=0$, the equation (\ref{qp2}) is 
       nothing but the classical second Painlev\'e equation $\PII$. 
       We call the equation (\ref{qp2}) \textsl{the quantum second Painlev\'e equation}.
       Also, putting $f_0=-(\varphi_0+\varphi_1)$, 
       $f_1=\varphi_0+\varphi_1+t-\psi^2$, $f_2=\psi$, (\ref{qp2}) 
       reduces to (\ref{0Nl=1}). However, in the same as $l=2n$ case, 
       it is not clear that there is some connection between the derivation 
       $\partial$ of $\mathcal{K}_1$ defined by the Hamiltonian $H_0$ in Section \ref{secqp}
       and the derivation $\partial$ of $\mathcal{F}_1$  defined by 
       the time evolution $T_1$ in this section.

\setcounter{section}{0}
\renewcommand\thesection{\Alph{section}}

\section{Appendix: Properties of the Hamiltonians $H_j$}

In the following, we introduce a family of Hamiltonians $H_1,\ldots,H_l$, by 
the diagram rotation. 
Namely,
\begin{equation}\label{eq:hamiltonians}
H_j:=\pi(H_{j-1}).
\end{equation}

In the classical case, these Hamiltonians have some remarkable 
properties in relation with the 
action of $W$. 

\begin{prop}\label{prop:hamiltonian1}
With respect to the action of the affine Weyl group $W$ 
(\ref{saction2}) and (\ref{saction}), 
the Hamiltonians have the following:

  (1) For $l=1,2n$, 
  \begin{equation}\label{eq:haw1}
    s_i(H_j)=H_j+\delta_{ij}k\frac{\alpha_j}{f_j}\quad (i,j=0,\ldots,l).
\end{equation}
  (2) For $l=2n+1$, 
  \begin{equation}\label{eq:haw2}
 s_i(H_j)=H_j+\delta_{ij}k\frac{\alpha_j}{f_j}x_j\quad (i,j=0,\ldots,l),
\end{equation}
where the index of $x_j$ (\ref{x0x1}) is regarded as in $\mathbb{Z}/2\mathbb{Z}$. 
In particular,  the Hamiltonian $H_j$ is invariant
with respect to the action of the Weyl group 
$W(A_l)=\langle s_0,\ldots,s_{j-1},s_{j+1},\ldots,s_l\rangle$.
\end{prop}

\begin{proof}
This proposition is a generalization of Proposition 
\ref{prop:hamiltonian0}.   
  From (\ref{haw1}), (\ref{haw2}), by using the definition
  of $H_j$ (\ref{eq:hamiltonians}) and the relation $\pi s_i=s_{i+1}\pi$, we obtain 
  the formulas (\ref{eq:haw1}), (\ref{eq:haw2}) respectively. \qed 
\end{proof}

\begin{prop}\label{prop:hamiltonian2}
(1) In the case of $A_{2n}^{(1)}$, for each $j=0,\ldots,2n$, 
one has 
\begin{equation}
H_{j+1}-H_j=k\sum_{r=1}^nf_{j+2r}-\frac{nk}{2n+1}x,
\end{equation}
where $x=f_0+f_1+\cdots+f_{2n}$ (\ref{x}).

(2) In the case of $A_{2n+1}^{(1)}$, for each $j=0,\ldots,2n+1$, 
one has
\begin{equation}
H_{j+1}-H_j=k\sum_{1\le r \le s \le n}^nf_{j+2r}f_{j+2s+1}-\frac{nk}{2n+1}\sum_{K\in S_2}f_K
+(-1)^j\frac{k}{4}\sum_{i=0}^l(-1)^i\alpha_i.
\end{equation}
\end{prop}

\begin{proof}
For $l=2,3$, we can prove through direct computations. 
For $l=2n$ ($l=2n,n\ge 2$), we have 
\begin{align*}
&H_0=\sum_{K\in S_3}f_K+\sum_{i=0}^{2n}\chi(\Gamma\backslash \{i\})f_i,
\\
&H_1=\sum_{K\in S_3}f_K+\sum_{i=0}^{2n}\pi(\chi(\Gamma\backslash \{i-1\}))f_i.
\end{align*}
Hence, 
\begin{equation}
H_1-H_0=\sum_{i=0}^{2n}\pi(\chi(\Gamma\backslash \{i-1\})
-\chi(\Gamma\backslash \{i\}))f_i.
\end{equation}
Computing $\pi(\chi(\Gamma\backslash \{i-1\})-\chi(\Gamma\backslash \{i\}))$ 
from the definition, we obtain 
\begin{equation}
\pi(\chi(\Gamma\backslash \{i-1\})-\chi(\Gamma\backslash \{i\}))=
\left\{
\begin{array}{ll}
\cfrac{-nk}{2n+1}& (i=0 \ \text{or}\  i=\text{odd})
\\[2mm]
\cfrac{(n+1)k}{2n+1}&(i\neq 0, i=\text{even})
\end{array}.
\right.
\end{equation}
Therefore, we obtain
\begin{equation}
H_1-H_0=\sum_{i=0}^{2n}\frac{-nk}{2n+1}f_i+\sum_{r=1}^nkf_{2r}=
k\sum_{r=1}^nf_{j+2r}-\frac{nk}{2n+1}x.
\end{equation}
Similarly we can show the formula in the case where $l=2n+1$ ($n\ge 2$). \qed
\end{proof}

In the classical case, we have $\tau$-functions $\tau_0,\ldots,\tau_l$ 
 such that $h_j=k(\log \tau_j)'$, 
where $h_j$ is the classical Hamiltonian corresponding to quantum Hamiltonian $H_j$. 
The affine Weyl group symmetry 
lifts to the level of $\tau$-functions. In fact,    
Proposition \ref{prop:hamiltonian1} and Proposition \ref{prop:hamiltonian2} 
illustrate how to lift the affine Weyl group symmetry to the level of 
$\tau$-functions. 

Unfortunately, the formulation of $\tau$-functions in the quantum case is not 
 completed yet and we hope to report on this in a near future.



\subsection*{Acknowledgements}

 The author is grateful to Koji Hasegawa 
 for suggestions and encouragements. He had drawn the author's attention to
  his preliminary draft
 written about a quantization of Kajiwara-Noumi-Yamada's realization \cite{kajiwarany} 
 of the affine 
 Weyl group of type $A_l^{(1)}$, and he suggested that a continuous limit of the 
 quantum discrete system constructed from that quantization would give the 
 quantum Painlev\'e system.  
 Also the author would like to thank Gen Kuroki for comments and 
 discussions. 
  

\end{document}